\documentclass[10pt.a4paper]{amsart}

\usepackage{amsmath}
\usepackage{amssymb}
\usepackage{amsfonts}

\usepackage{graphicx}
\usepackage{epstopdf}
\numberwithin{equation}{section}
\newtheorem{theorem}{Theorem}[section]

\newtheorem{proposition}[theorem]{Proposition}

\theoremstyle{remark}
\newtheorem{example}{Example}[section]

\theoremstyle{definition}
\newtheorem{definition}{Definition}[section]

\theoremstyle{remark}

\theoremstyle{remark}
\newtheorem{remark}{Remark}[section]

\begin{document}

\newcommand{\spt}{\,\mathrm{supp}\,}
\newcommand{\xxi}{\langle\xi\rangle}
\newcommand{\xx}{\langle x\rangle}
\newcommand{\yy}{\langle y\rangle}
\newcommand{\dint}{\int\!\!\int}
\newcommand{\triple}[1]{{|\!|\!|#1|\!|\!|}}
\newcommand{\Sph}{{\mathbb S}}
\newcommand{\LII}[3]{L_{#1}^{#2}L_{\phantom{I}}^{#3}}
\newcommand{\LI}[2]{\LII{I}{#1}{#2}}

\title[Schr\"odinger with potential]%
{Some remarks on the Schr\"odinger equation\\
with a potential in $L^{r}_{t}L^{s}_{x}$}

\date{\today}      

\thanks{{\em Support.} The authors are partially supported by
the Research Training Network (RTN) HYKE
and by grant
HPRN-CT-2002-00282  from the European Union. 
The third author is supported 
also by INDAM}

\subjclass[2000]{
35B40, 35B25, 35B65, 35Q40, 35Q55.}

\keywords{%
Schr\"odinger equation,
time-dependent potential,
Strichartz estimates.
}
\maketitle


\centerline{Piero D'Ancona}
\vspace{0.2cm}
\centerline{\em Universit\`a di Roma ``La Sapienza''}
\centerline{\em Dipartimento di Matematica,}
\centerline{\em Piazzale A.~Moro 2, I-00185 Roma, Italy,}
\centerline{\em email: dancona@mat.uniroma1.it}
\centerline{tel/fax: +39-06-4991-3092}
\vspace{0.4cm}
\centerline{Vittoria Pierfelice and Nicola Visciglia}
\vspace{0.2cm}
\centerline{\em Universit\`a di Pisa}
\centerline{\em Dipartimento di Matematica}
\centerline{\em Via F. Buonarroti 2, I-56127 Pisa, Italy}
\centerline{\em email: pierfelice@dm.unipi.it, viscigli@mail.dm.unipi.it}

\vspace{0.2cm}

{\bf Abstract.}
{\em We study the dispersive properties of the linear Schr\"odinger equation
with a time-dependent potential $V(t,x)$. We show that  an 
appropriate integrability
condition in space and time on $V$, i.e. the boundedness of a 
suitable $L^{r}_{t}L^{s}_{x}$ norm,
is sufficient to prove the full set of Strichartz estimates.
We also construct several counterexamples which show that our 
assumptions are optimal,
both for local and for global Strichartz estimates, in the class 
of large unsigned 
potentials $V\in L^{r}_tL^{s}_x$.}

\section{Introduction}\label{sec.introd}

Dispersive properties of evolution equations have become in recent
years a crucial tool in the study of a variety of questions, including 
local and
global existence for nonlinear equations, well posedness in Sobolev 
spaces of low order,
scattering theory and many others. In particular, the free 
Schr\"odinger equation
\begin{equation}\label{eq.freesch}
    i\partial_{t}u-\Delta u=0
\end{equation}
exhibits a rich set of dispersive and smoothing properties. 
The basic ones are the
\emph{$L^{p}-L^{q}$ estimates}
\begin{equation}\label{eq.disp}
    \|u(t,\cdot)\|_{L^{q}}\leq \frac C {t^{\sigma}}\|u(0,\cdot)\|_{L^{p}},
    \qquad
    \frac 1 p +\frac 1 q=1,\quad p\in[1,2],\quad \sigma=\frac n 2-\frac n q,
\end{equation}
and the \emph{Strichartz estimates}, which we shall recall 
shortly (see Definition \ref{def.strich} below). 

The usual route to prove
these estimates is: first deduce the $L^{p}-L^{q}$ estimates 
from the explicit expression of the fundamental solution; then, 
use them to derive the Strichartz estimates via suitable 
functional analytic arguments. Standard references on 
the subject are \cite{GinibreVelo95-generstric} and 
\cite{KeelTao98-endpoinstric}; see also \cite{Cazenave03-SLS}. 
However, in more general situations this approach is not feasible, and indeed, 
for some generalizations of the Schr\"odinger equation, 
Strichartz estimates may hold even when pointwise 
estimates fail (see e.g. \cite{Bourgain93-NLS}, 
\cite{Bourgain93-NLSfourier},  \cite{BurqGerardTzvetkov04-strichcompact}). 
Thus it appears that Strichartz estimates have a 
more fundamental nature and a greater generality.

Here we shall focus on the Schr\"odinger equation 
with a time-dependent potential
\begin{equation}\label{eq.schrpot}
    i\partial_{t}u-\Delta u+V(t,x)u=0
\end{equation}
and its inhomogeneous version with a source term.
The great interest of this equation, both from the physical 
and from the mathematical point of view, is well known. 
Many results on the dispersive properties are available 
when the potential $V=V(x)$ depends only on space 
variables (see, among others, 
\cite{GoldbergSchlag}, \cite{JourneSofferSogge91-decayschroed},
\cite{RodnianskiSchlag04-disp} and the references therein). 
On the other hand, the time-dependent case is much more 
delicate; almost all available results are of a perturbative 
nature, requiring some smallness of the potential $V(t,x)$ 
(see \cite{GoldbergSchlag}, \cite{RodnianskiSchlag04-disp},
and, for small potentials of very low regularity,
\cite{GI}, \cite{Pierfelice}; concerning the case of 
time-periodic potentials, see
\cite{H},  and \cite{Yajima84-timeperiodic};
see also \cite{costin}).

Our goal here is to show that, by purely elementary arguments 
based on integrability properties of the potential (as opposed 
to the "global" smallness required in the above mentioned 
results), it is possible to obtain a great deal of information 
on the behaviour of the solution, and to prove the Strichartz 
estimates for a wide class of large potentials with no definite 
sign. Of course, the usual obstructions are present also in this 
general situation: existence of standing waves, rescaling and 
pseduoconformal symmetry of the equation. Using these, we 
are able to show that our conditions are also necessary, at 
least in the class of potentials under consideration.

Let us recall the classical Strichartz estimates for the Schr\"odinger 
equation, and introduce some notations. We use a prime to denote
conjugate indices; moreover, for any subinterval $I$ of $\mathbb{R}$ 
(bounded or unbounded) we define the mixed space-time norms
\begin{equation}\label{eq.spt}
    \|u\|_{\LI p q}\equiv \left(\int_{I}\|u(t,\cdot)\|_{L^{q}
    (\mathbb{R}^{n})}^{p}dt \right)^{1/p}
    \end{equation}
and when $I=[0,+\infty[$ we write simply $L^{p}L^{q}$ i
n place of $\LI p q$. Similarly, we shall write
\begin{equation}\label{eq.CL2}
    C_{I}L^{p}\equiv C(I;L^{p}),\qquad CL^{p}\equiv C([0,+\infty[;L^{p})
\end{equation}
for $1\leq p\leq\infty$.

\begin{definition}\label{def.strich}
Let $n\geq2$. The pair $(p,q)$ is said to be \emph{(Schr\"odinger) 
admissible} if
\begin{equation}\label{eq.strich1}
    \frac 1 p + \frac n {2q} = \frac n 4,\qquad p,q\in[2,\infty],\qquad
    (n,p,q)\neq(2,2,\infty).
\end{equation}
The Strichartz estimates can be stated as follows:
for all admissible couples $(p,q)$ and
$(\tilde p,\tilde q)$ there exists a constant $C(p,\tilde p)$ 
such that, for all interval $I\subseteq\mathbb{R}$ (bounded 
or unbounded), for all functions $u_{0}(x)\in L^{2}(\mathbb{R}^{n})$,
and $F(t,x)\in \LI {\tilde p'}{\tilde q'}$ the following inequalities hold:
\begin{equation}\label{eq.strich}
    \biggl\| e^{it\Delta}u_{0}\;\biggr\|_{\LI p q}\leq C(p,\tilde p)\;\|u_{0}\|_{L^{2}}
\end{equation}
\begin{equation}\label{eq.strichF}
    \biggl\|\int_{0}^{t}e^{i(t-s)\Delta}F(s)\;ds\biggr\|_{\LI p q}\leq 
       C(p,\tilde p)  \;  \bigl\| F\bigr\|_{\LI{\tilde p'}{\tilde q'}}
\end{equation}
Note that the constant is independent of the interval $I$. 

Clearly, when $n\geq3$ the constant
can be taken also independent of $p$ and $\tilde p$: 
we shall denote this universal constant (which depends 
now only on the space dimension $n$) by $C_{0}$. 
When $n=2$, the constant is unbounded as $p\downarrow 2$ 
or $\tilde p\downarrow 2$.
\end{definition}

Here $e^{it\Delta}$ is the unitary operator
\begin{equation}\label{eq.defschr}
    e^{it\Delta}f= \int_{\mathbb{R}^{n}}
        \frac{e^{-\frac{i|x-y|^{2}}{4t}}}{(4\pi i t)^{n/2}}\;f(y)\;dy,\qquad
    \int_{0}^{t}e^{i(t-s)\Delta}F(s)\;ds= 
       \int_{0}^{t}\int_{\mathbb R^n}
        \frac{e^{-\frac{i|x-y|^{2}}{4(t-s)}}}{(4\pi i (t-s))^{n/2}}\;F(s,y)\;dy\;ds
\end{equation}
which is properly defined on $L^{2}$ but can be extended to 

different $L^{p}$ spaces
using e.g.\ these explicit expressions.

Consider now the differential equation
\begin{equation}\label{eq.schrV}
        i\partial_{t}u-\Delta u+V(t,x) u=F(t,x) , \qquad
        u(0,x)=u_{0}(x).
\end{equation}
For low regularity solutions, it is customary to replace 
\eqref{eq.schrV} with the integral equation
\begin{equation}\label{eq.schrI}
    u(t,x)=e^{it\Delta }u_{0}(x)+\int_{0}^{t}e^{i(t-s)\Delta}[F(s)-V(s)u(s)]\;ds.
\end{equation}
The two formulations are equivalent under very mild 
assumptions on the class of solutions; we shall not 
discuss this problem here, instead we shall use the 
integral formulation exclusively.

We can now state our first result:

\begin{theorem}\label{th.1}
Let $n\geq2$, 
let $I$ be either the interval $[0,T]$ or $[0,+\infty[$, 
and assume $V(t,x)$ is
a real valued potential belonging to
\begin{equation}\label{eq.assV}
    V(t,x)\in \LI r s,\qquad \frac 1 r + \frac n {2s}=1
\end{equation}
for some fixed $r\in[1,\infty[$ and $s\in ]n/2,\infty]$. 
Let $u_{0}\in L^{2}$ and 
$F\in\LI{\tilde p'}{\tilde q'}$ for some admissible pair 
$({\tilde p'},{\tilde q'})$.

Then the integral equation \eqref{eq.schrI} has
a unique solution $u\in C_{I}L^{2}$ which belongs to 
$\LI p q$ for all admissible pairs $(p,q)$
and satisfies the Strichartz estimates
\begin{equation}\label{eq.strV}
    \|u\|_{\LI p q} \leq C_{V} \; \|u_{0}\|_{L^{2}}+C_{V}\|F\|_{\LI {\tilde p'}{\tilde q'}}.
\end{equation}
When $n\geq3$, the constant $C_{V}$ can be estimated 
by $k(1+2C_{0})^{k}$, where $C_{0}$ is the Strichartz 
constant for the free equation, while $k$ is an integer such that
the interval $I$ can be partitioned in $k$ subintervals $J$ 
with the property
$\|V\|_{\LII{J} r s}\leq (2C_{0})^{-1}$. A similar statement 
holds when $n=2$, provided we
replace $C_{0}$ by $C({p,\tilde p)}$.

Finally, when $F\equiv0$ the solution satisfies the conservation of energy
\begin{equation}\label{eq.cons}
    \|u(t)\|_{L^{2}}\equiv\|u_{0}\|_{L^{2}},\qquad t\in I.
\end{equation}
\end{theorem}

\begin{remark}\label{rem.emph}
We emphasize that the potentials $V(t,x)$ considered in 
Theorem \ref{th.1} may be both large and change sign. 
The usual smallness assumption is replaced here by the 
integrability condition with respect to time, which ensures 
smallness of $V$ on sufficiently small time intervals, and for $t>>1$.
\end{remark}

\begin{remark}\label{rem.sum}
By iterating the argument of the proof, it is easy to extend 
the above result to any potential of the form
\begin{equation*}
    V=V_{1}+\cdots+V_{k}
\end{equation*}
where $V_{1},\dots,V_{k}$ satisfy assumption \eqref{eq.assV}, 
with possibly different indices $r_{j}, s_{j}$. 
\end{remark}

\begin{remark}\label{rem.ineq}
Note that when $I$ is a bounded interval, assumption 
\eqref{eq.assV} can be relaxed to
\begin{equation}\label{eq.assV2}
    V(t,x)\in \LI r s,\qquad \frac 1 r + \frac n {2s}\leq 1;
\end{equation}
indeed, from \eqref{eq.assV2}, using H\"older's inequality 
in time we can easily show that also \eqref{eq.assV} holds, 
for a smaller value of $r$ and the same value of $s$.

Thus in particular we see that the existence part of our 
theorem extends a result of Yajima 
\cite{Yajima87}, who proved that the equation \eqref{eq.schrI} 
(or \eqref{eq.schrV}) is well posed in $L^{2}$ with conservation 
of energy, provided the potential $V$ satisfies
\begin{equation}\label{eq.yaj}
    V=V_{1}+V_{2},\qquad V_{1}\in\LI  r s,\quad V_{2}\in\LI \infty \beta
\end{equation}
with $\beta>1$ and
\begin{equation}\label{eq.yaj2}
    \frac 1 r+\frac n {2s}<1
\end{equation}
(see also the preceding remark).
\end{remark}

When the potential $V(t,x)$ belongs to $\LI \infty {n/2}$, i.e., 
in the limiting case
of Theorem \ref{th.1}, the result can not be true; indeed, this 
case includes
the static potentials $V(x)\in L^{n/2}$ without any positivity or smallness
assumption. We mention that even for a nonnegative potential in 
$ L^{n/2}$ it is not known if  the Strichartz
estimates are valid in general. The best result in this direction is 
due to Rodnianski and Schlag
\cite{RodnianskiSchlag04-disp} who considered bounded 
potentials defined on $\mathbb{R}^n$
satisfying the estimate 
$|V(x)|\leq C| x|^{-2-\epsilon}$ for $|x|$ large enough. 
However, in the limiting case we can prove a partial 
substitute of Theorem \ref{th.1}.
To simplify our statement we introduce the following definition:

\begin{definition}\label{def.V}
Let $V(x)$ be a real valued function such that
$$H=\Delta-V(x)$$
has a selfadjoint extension.
We say that the potential $V(x)$ is  of \emph{Strichartz type} if 
for all bounded time interval $I=[0,T]$, for all $u_{0}\in L^{2}$
and $F\in\LI {\tilde p'}{\tilde q'}$ with $(\tilde p,\tilde q)$ admissible,
the integral equation
\begin{equation}\label{eq.schrIVx}
    u(t,x)=e^{itH}u_{0}+\int_{0}^{t}e^{i(t-s)H}F(s)\;ds
\end{equation}
has a unique solution $u(t,x)\in C_{I}L^{2}$  satisfying the
Strichartz estimates
\begin{equation}\label{eq.strVx}
    \|u\|_{\LI a b} \leq C(I,V) \; \|u_{0}\|_{L^{2}}+C(I,V)\;\|F\|_{\LI {\tilde p'}{\tilde q'}}
\end{equation}
for all admissible pairs $(a,b)$.
\end{definition}

Then we have:

\begin{theorem}\label{th.2}
Let $n\geq3$, let $I$ be a bounded interval $[0,T]$ and let 
$V(t,x)\in C_{I}L^{n/2}$. Assume that
for each $t\in I$, $V(t,\cdot)$ is of Strichartz type, while the 
functions $u_{0}$ and $F(t,x)$ are as in Theorem \ref{th.1}. 
Then the conclusion of Theorem \ref{th.1} holds 
true (local Strichartz estimates). 

The result holds also in the case $I=[0,\infty[$ (global Strichartz 
estimates) under the additional assumption: there exists 
$T_{0}>0$ such that $\|V(t,\cdot)\|_{L^{n/2}}\leq (2C_{0})^{-1}$ for 
$t>T_{0}$.
\end{theorem}

\begin{remark}\label{rem.nbd}
By simple modifications in the proof, Theorem \ref{th.2} can 
be extended to any potential of the form
$$V(t,x)=V_{1}(t,x)+V_{2}(t,x),$$
with $V_{1}$ as in the theorem while $V_{2}\in \LI \infty {n/2}$ satisfies
$$\|V_{2}\|_{\LI \infty {n/2}}\leq \varepsilon(V_{1})$$
for a suitable small constant $\epsilon(V_{1})$ depending only on $V_{1}$.
\end{remark}

\begin{example}\label{ex.NLS}
To illustrate a possible use of Theorem \ref{th.1}, consider 
the semilinear Schr\"odinger equation
\begin{equation}\label{eq.NLS}
    i\partial_{t}u-\Delta u+f(u)u=0,\qquad |f(u)|\leq 
    C|u|^{\gamma},\quad \gamma>1,
\end{equation}
$f$ real valued,
which includes both focusing and defocusing equations with 
a power nonlinearity. Then we may regard \eqref{eq.NLS} 
as a Schr\"odinger equation with a time dependent potential
\begin{equation*}\label{eq.VNLS}
    V(t,x)={f(u(t,x))}.
\end{equation*}
We see that $V$ satisfies the assumptions of Theorem \ref{th.1} provided
\begin{equation}\label{eq.assu}
    u\in L^{a}L^{b},\qquad \frac 1 a+\frac n {2b}=
    \frac 1 {\gamma}, \quad a<\infty.
\end{equation}
Thus any solution satisfying \eqref{eq.assu} satsfies the full 
set of Strichartz estimates. 

For instance, in the case of the (focusing or defocusing) quintic 
Schr\"odinger equation in three dimensions, any solution
$u\in L^{10}L^{10}$ satisfies the Strichartz estimates; this 
was the first step in the proof of the global well posedness 
for the radial defocusing three dimensional quintic in 
\cite{Bourgain99-defoc}.
\end{example}

\begin{example}\label{ex.goldberg}
We give a simple application of Theorem \ref{th.2}. Consider 
a real valued potential 
$V\in CL^\frac n{2}$ and assume it satisfies the bounds
\begin{equation}\label{eq.VR3}
    0\leq V(t,x)\leq\frac C {(1+|x|)^{2+\delta}},\qquad 
    x\in \mathbb{R}^{n},\quad n\geq3
\end{equation}
for some $C,\delta>0$. Then we can prove that $V(t,x)$ 
satisfies the assumptions of Theorem \ref{th.2} and hence the
local Strichartz estimates hold (and also the global ones, 
under the additional  assumption of smallness at infinity). 

Indeed, let $W(x)=V(t_{0},x)$ for an arbitrary fixed $t_{0}$; 
we must show that $W(x)$ is of Strichartz type. The existence 
part of the definition is trivial; let us prove the estimates. 
Consider the operator $H=-\Delta+W(x)$; $H$ has a 
unique selfadjoint extension by standard results, with 
spectrum contained in $[0, +\infty[$; by Theorem XIII.58 
in \cite{ReedSimonIV} $H$ has no strictly positive eigenvalues, 
since $W$ is bounded and decays as $|x|^{-2-\delta}$ at 
infinity; 0 is certainly not an eigenvalue since $Hf=0$ implies 
$f=0$ easily. This implies that the operator $H$ has a purely 
continuous spectrum. Now Theorem 1.4 in 
\cite{RodnianskiSchlag04-disp} states that $P_{c}e^{itH}$ 
satisfies the full set of Strichartz estimates when the potential 
is bounded and decays faster than $|x|^{-2}$ at infinity; 
here $P_{c}$ is the projection on the continuous subspace 
of $L^{2}$ for $H$, which coincides with all of $L^{2}$ as 
we have just proved. In conclusion, $W(x)=V(t_{0},x)$ is of Strichartz type as 
claimed.
\end{example}

\begin{remark}\label{rem.resc}
Condition \eqref{eq.assV} is quite natural, in view of the 
following argument: the standard rescaling 
$u_{\epsilon}(t,x)=u(\epsilon^{2}t,\epsilon x)$ 
takes equation \eqref{eq.schrpot} into the equation
\begin{equation}\label{eq.rescaled}
        i\partial_{t}u_{\epsilon}-\Delta u_{\epsilon}
	+V_{\epsilon}(t,x)u_{\epsilon}=0,\qquad
        V_{\epsilon}(t,x)=\epsilon^{2}V(\epsilon^{2}t,\epsilon x),
\end{equation}
and we have
\begin{equation}
    \|V_{\epsilon}\|_{L^{r}L^{s}}=
    \epsilon^{2\left(1-\frac {1}{r}-\frac {n} {2s}\right)}\|V\|_{L^{r}L^{s}}
\end{equation}
so that the $L^{r}L^{s}$ norm of $V_{\epsilon}$ is independent 
of $\epsilon$ precisely when $r,s$ satisfy \eqref{eq.assV}. 

Indeed, by a suitable use of rescaling arguments, it is possible 
to show that the condition 
$1/r+n/(2s)=1$ is \emph{necessary} in order that the global 
Strichartz estimates be true for any potential belonging to the 
classes $L^{r}L^{s}$ (see Theorem \ref{th.3}  below).

Concerning the \emph{local} Strichartz estimates, the situation 
is more interesting. When 
$1/r+n/(2s)<1$, as already observed in Remark \ref{rem.ineq}, 
the local Strichartz estimates are an elementary consequence of 
Theorem \ref{th.1}.  On the other hand, when $1/r+n/(2s)>1$,
it is possible to show that the \emph{local} Strichartz estimates fail. 
This case is more delicate; actually it is not even clear if equation 
\eqref{eq.schrpot} is well posed in $L^{2}$ under this assumption on $V$.
\end{remark}

We collect our counterexamples in the following theorem, concerning
the homogeneous equation
\begin{equation}\label{eq.schVhom}
    iu_{t}-\Delta u+V(t,x)u=0.
\end{equation}
Note that the case $(r,s)=(\infty,n/2)$ is almost trivial since it 
is based on the construction of a standing wave for 
\eqref{eq.schVhom}; we state it in some length both 
for completeness, and because the remaining counterexamples 
are based on it. Thus,  in the proof of Theorem \ref{th.3} 
it is essential to use potentials which change sign. 

\begin{theorem}\label{th.3}
Let $n\geq2$. Then we have the following.

(i) (Case $r=\infty$) We can construct a potential 
$W(x) \in C^{\infty}_{0}(\mathbb{R}^{n})$
and a function $u_{0}\in H^{s}$ for all $s>0$ such that
\begin{equation}\label{eq.inf}
    -\Delta u_{0}+W(x)u_{0}+u_{0}=0.
\end{equation}
Hence the function $u(t,x)=e^{-it}u_{0}(x)\in CL^{2}$ solves 
\eqref{eq.schVhom} with
$$V(t,x)\equiv W(x)\in L^{\infty}([0,+\infty[;L^{n/2}(\mathbb{R}^{n})),$$
and does not belong to the space 
$L^{p}([0,+\infty[;L^{q})$ for all admissible pairs
$(p,q)\neq(\infty,2)$. In other words, there exists a potential 
$V(t,x)$ belonging to $ L^{\infty}L^{s}$ for all $s\in[1,\infty]$ 
such that the global Strichartz estimates \eqref{eq.strV} on 
$I=[0,+\infty[$ do not hold for equation \eqref{eq.schVhom}.

(ii) (Counterexamples to global estimates) 
For every pair $(r,s)$ with $r\in[1,\infty[$, $s\in ]n/2,\infty]$ and
\begin{equation}\label{eq.contro}
     \frac 1 r + \frac n {2s}\neq1,
\end{equation}
we can construct a potential $V(t,x)\in L^{r}([0,+\infty[;L^{s})$ 
and a sequence of solutions $u_{k}(t,x)\in C([0,+\infty[{};L^{2})$ 
to equation \eqref{eq.schVhom} such that
\begin{equation}\label{eq.viol}
    \lim_{k\to\infty}\frac{\|u_{k}\|_{L^{p}L^{q}}}{\|u_{k}(0)\|_{L^{2}}}=
    \infty\quad
    \hbox{for every admissible pair}\quad (p,q)\neq(\infty,2).
\end{equation}

(iii) (Counterexamples to local estimates) 
For every pair $(r,s)$ with $r\in[1,\infty[$, $s\in ]n/2,\infty]$ and
\begin{equation}\label{eq.contro2}
     \frac 1 r + \frac n {2s}>1,
\end{equation}
we can construct, on any given bounded time interval $I=[0,T]$, 
a potential 
$V(t,x)\in L^{r}([0,T];L^{s})$ and a sequence of solutions 
$u_{k}(t,x)\in C([0,T];L^{2})$ to equation \eqref{eq.schVhom} such that
\begin{equation}\label{eq.viol2}
    \lim_{k\to\infty}\frac{\|u_{k}\|_{\LI p q}}{\|u_{k}(0)\|_{L^{2}}}=\infty\quad
    \hbox{for every admissible pair}\quad (p,q)\neq(\infty,2).
\end{equation}
\end{theorem}

We conclude the paper with a result showing that, at least 
for a restricted range of indices $r,s$, the conclusion of 
Theorem \ref{th.3}, part (iii), can be improved in an essential way. 
While the above theorem was based on suitable rescaling 
arguments, Proposition \ref{prop.4} exploits the 
\emph{pseudoconformal} invariance of the Schr\"odinger equation.

\begin{proposition}\label{prop.4}
Let $n\geq2$, and assume $r\in[1,\infty[$ and $s\in]n/2,n[$ satisfy
\begin{equation}\label{eq.condcontr}
    \frac 1{2r}+\frac n {2s}>1.
\end{equation}
Then we can construct a potential 
$V(t,x)\in L^{r}(0,1;L^{s}(\mathbb{R}^{n}))$ and a solution 
$u(t,x)\in C([0,1];L^{2})$ to equation \eqref{eq.schVhom}
such that, for all admissible pairs $(p,q)$ with $p<\infty$, 
and for any $0<T<1$, we have
\begin{equation*}
    u\in L^{p}(0,T;L^{q}(\mathbb{R}^{n}))\qquad\hbox{but}\qquad
    u\not\in{L^{p}(0,1;L^{q}(\mathbb{R}^{n}))}.
\end{equation*}
\end{proposition}

\section{Proof of Theorem \ref{th.1} }\label{sec.th1}

We shall consider in detail only the case $n\geq3$; in the case $n=2$, when the
endpoint fails, it is sufficient to replace in the following arguments the space
$ \LII J 2 {\frac{2n}{n-2}}$ with any $\LII J p q$ with $q$ arbitrarily large.

We distinguish two cases, according to the value of $r\in[1,\infty[$. 

\subsection{Case A: $r\in[2,\infty[$}\label{ssec.caseA}

Consider a small interval $J=[0,\delta]$, and let $Z$ be the Banach space
$$Z=C_{J}L^{2}\cap \LII J 2 {\frac{2n}{n-2}},\qquad
     \|v\|_{Z}:=\max \left\{ \|v\|_{\LII J \infty 2},  \;\;  \|v\|_{\LII J 2 
     {\frac{2n}{n-2}}} \right\}.$$
Notice that, by interpolation, $Z$ is embedded in all admissible 
spaces $\LII J p q$.

For any $v(t,x)\in Z$ we define the mapping
\begin{equation}\label{eq.phi}
    \Phi(v)=e^{it\Delta }u_{0}+\int_{0}^{t}e^{i(t-s)\Delta}[F(s)-V(s)v(s)]\;ds.
\end{equation}
A direct application of \eqref{eq.strich}, \eqref{eq.strichF} gives
\begin{equation}\label{eq.str1}
    \|\Phi(v)\|_{\LII J p q}\leq C_{0}\|u_{0}\|_{L^{2}}+
       C_{0}\|V v\|_{\LII J {p_{0}'} {q_{0}'}}+
       C_{0}\|F\|_{\LII J {\tilde p'}{\tilde q'}}
\end{equation}
for all admissible $(p,q)$, $(p_{0},q_{0})$, $(\tilde p, \tilde q)$, 
and by H\"older's inequality
we can write
\begin{equation}\label{eq.hold}
    \|\Phi(v)\|_{\LII J p q}\leq C_{0}\|u_{0}\|_{L^{2}}+
       C_{0}\|V\|_{\LII J r s}  \|v\|_{\LII J {2} {\frac{2n}{n-2}}}+
       C_{0}\|F\|_{\LII J {\tilde p'}{\tilde q'}}
\end{equation}
provided we choose $p_{0}, q_{0}$ such that
$$\frac 1 {p_{0}}=\frac 1 2 -\frac 1 r,\qquad
     \frac 1 {q_{0}}=\frac {n+2} {2n} -\frac 1 s.$$
Note that
$$\frac 1 {p_{0}}+\frac n {2q_{0}}=
    \frac 1 2+\frac {n+2} 4-\left(\frac 1 r +\frac n {2s}\right)
       \equiv
    \frac 1 2+\frac {n+2} 4-1
        \equiv
      \frac n 4$$
by our assumptions on $r,s$, and moreover
$$r\in[2,\infty[\quad \implies \quad p_{0}\in[2,\infty[$$
so that our choice of $p_{0},q_{0}$ always gives an admissible pair
in the case under consideration.

In particular, choosing $(p,q)=(\infty,2)$ or $(2,2n/(n-2))$, we obtain
\begin{equation}\label{eq.hold2}
    \|\Phi(v)\|_{Z}\leq C_{0}\|u_{0}\|_{L^{2}}+
       C_{0}\|V\|_{\LII J r s}  \|v\|_{Z} +
       C_{0}\|F\|_{\LII J {\tilde p'}{\tilde q'}}
\end{equation}

Thus $\Phi(v)$ belongs to all the admissible spaces $\LII J p q$, 
and to prove that
$\Phi(v)$ belongs to $Z$ it remains only to show that $u$ 
is continuous with
values in $L^{2}$. But this is an immediate consequence of the 
following simple remark:

\begin{remark}\label{rem.CL2}
Let $G(t,x)\in\LII J {a'} {b'}$ with $(a,b)$ admissible.
Then the function
$$w(t,x)=\int_{0}^{t}e^{i(t-s)\Delta}G(s)\;ds$$
belongs to $C_{J}L^{2}$. Indeed, this is certainly true if we 
know in addition that $G$ is
a smooth function, compactly supported in $x$ for each $t$. 
If we approximate $G$ by a
sequence of such functions  $G_{j}$ in the $\LII J {a'} {b'}$ norm, 
the Strichartz estimates imply that
the corresponding functions $w_{j}$ converge in $L^{\infty}L^{2}$, 
whence the claim follows.
\end{remark}

We have thus constructed a mapping $\Phi:Z\to Z$. 
Assume now the length $\delta$ of the
interval $J$ is chosen so small that\
\begin{equation}\label{eq.smallV}
    C_{0}\|V\|_{\LII J r s}\leq \frac1 2;
\end{equation}
this is certainly possible since $r<\infty$. With this choice we 
obtain immediately
two consequences: first of all, the mapping $\Phi$ is a 
contraction on $Z$ and hence has a unique fixed point 
$v(t,x)$ which is the required solution; second, $v$ satisfies
\begin{equation}\label{eq.contr}
    \|v\|_{\LII J p q}\leq C_{0}\|u_{0}\|_{L^{2}}+
       \frac12\|v\|_{\LII J {p} {q}}+
       C_{0}\|F\|_{\LII J {\tilde p'}{\tilde q'}}
\end{equation}
whence we obtain
\begin{equation}\label{eq.smallstr}
        \|v\|_{\LII J p q}\leq 2C_{0}\|u_{0}\|_{L^{2}}+
               2C_{0}\|F\|_{\LII J {\tilde p'}{\tilde q'}}
\end{equation}

It is clear that the above argument applies on any 
subinterval $J=[t_{0},t_{1}]\subseteq I$ on which
a condition like \eqref{eq.smallV} holds; of course, we will 
obtain an estimate of the form
\begin{equation}\label{eq.smallstrj}
        \|v\|_{\LII J p q}\leq 2C_{0}\|v(t_{0})\|_{L^{2}}+
               2C_{0}\|F\|_{\LII J {\tilde p'}{\tilde q'}}.
\end{equation}
Notice also that \eqref{eq.smallstrj} implies in particular
\begin{equation}\label{eq.induct}
        \|v(t_{1})\|_{L^{2}}\leq 2C_{0}\|v(t_{0})\|_{L^{2}}+
               2C_{0}\|F\|_{\LII J {\tilde p'}{\tilde q'}}.
\end{equation}

Now we can partition the interval $I$ (bounded or unbounded) 
in a finite number of subintervals
on which condition \eqref{eq.smallV} holds. Applying inductively 
the estimates \eqref{eq.smallstrj}
and \eqref{eq.induct} we easily obtain \eqref{eq.strV} and the 
claimed estimate for the Strichartz constant.

The last remark \eqref{eq.cons} concerning the conservation of 
energy can be proved by approximation as follows: let 
$V_{j}(t,x)$ be a sequence of real valued smooth potentials, 
compactly
supported in $x$, and let $v_{j}$ be the corresponding 
solutions; then the differences
$w_{j}=v-v_{j}$ satisfy (in suitable integral sense)
$$i\partial_{t}w_{j}-\Delta w_{j}+Vw_{j}=(V-V_{j})v_{j}\equiv F_{j}.$$
Now we observe that the smooth solutions $v_{j}$ have 
a conserved energy;
moreover, we can choose the approximating potentials 
$V_{j}$ in such a way that
they converge to $V$ in $\LI r s$ and their Strichartz 
constants do not exceed the above constructed 
constant for $V$. Indeed, if we can partition $I$ in a 
finite set of subintervals
satisfying \eqref{eq.smallV}, we can choose exactly the 
same subintervals for each $V_{j}$
provided we construct $V_{j}$ by a convolution with standard 
mollifiers, so that their Lebesgue
norm does not increase. In conclusion, the $v_{j}$ satisfy 
uniform Strichartz estimates, and this
implies that the nonhomogeneous terms $F_{j}=(V-V_{j})v_{j}$ 
tend to 0 in the (dual) admissible spaces, by estimates 
identical to the above ones. Thus in particular $w_{j}\to 0$ 
in $L^{\infty}L^{2}$ and this shows that also $v(t,x)$ 
satisfies the conservation of energy.

\subsection{Case B: $r\in[1,2]$}\label{ssec.caseB}

The method in this case is quite similar to the above one, 
but instead of \eqref{eq.str1}
we use the estimate
\begin{equation}\label{eq.str2}
    \|\Phi(v)\|_{\LII J p q}\leq C_{0}\|u_{0}\|_{L^{2}}+
       C_{0}\|V v\|_{\LII J {r} {\frac{2s}{s+2}}}+
       C_{0}\|F\|_{\LII J {\tilde p'}{\tilde q'}}
\end{equation}
where $(p,q)$ and $({\tilde p},{\tilde q})$ are arbitrary 
admissible pairs, while
the pair $(r,2s/(s+2))$ is the dual of $(r',2s/(s-2))$ and 
this last pair is admissible since
$$\frac 1{r'}+\frac n 2\cdot\frac {s-2}{2s}=
        \frac n {2s}+\frac n 2\cdot\frac{s-2}{2s}=\frac n 4$$
where  we have used the assumption $1/r+n/(2s)=1$; 
notice also that $r\in[1,2]$ and
hence $2s/(s+2)\in [1,2]$ too.

Thus by H\"older's inequality we obtain
\begin{equation}\label{eq.str3}
    \|\Phi(v)\|_{\LII J p q}\leq C_{0}\|u_{0}\|_{L^{2}}+
       C_{0}\|V\|_{\LII J r s} \|v\|_{\LII J {\infty} {2}}+
       C_{0}\|F\|_{\LII J {\tilde p'}{\tilde q'}}
\end{equation}
and choosing $(p,q)=(\infty,2)$ or $(2,2n/(n-2))$  and proceeding as above
we arrive at
\begin{equation}\label{eq.str4}
        \|\Phi(v)\|_{Z}\leq C_{0}\|u_{0}\|_{L^{2}}+
       \frac 1 2\|v\|_{Z}+
       C_{0}\|F\|_{\LII J {\tilde p'}{\tilde q'}}.
\end{equation}
From this point on, the proof is identical to the first case.

\section{Proof of Theorem \ref{th.2} }\label{sec.th2}

The proof follows the same lines as the preceding one; 
indeed, the continuity
in time of the potential allows to consider $V(t,x)$ as a small perturbation of
$V(t_{0},x)$ for $t$ near $t_{0}$. 

Let $J=[0,\delta]$ be a small interval, and consider again the space
$$Z=C_{J}L^{2}\cap \LII J 2 {\frac{2n}{n-2}},\qquad
     \|v\|_{Z}:=\max \left\{ \|v\|_{\LII J \infty 2},  \;\;  \|v\|_{\LII J 2 
     {\frac{2n}{n-2}}} \right\}.$$
On $Z$ we construct a map $\Phi$ defined as follows:
\begin{equation}\label{eq.phi2}
    \Phi(v)=e^{itH}u_{0}
       +\int_{0}^{t}e^{i(t-s)H}[F(s)-W(s)v(s)]\;ds,
\end{equation}
where
\begin{equation}\label{eq.H}
    H=\Delta-V(0,x),\qquad W(t,x)=V(t,x)-V(0,x).
\end{equation}
We have used the assumption that $V(0,x)$ is of Strichartz 
type (Definition \ref{def.V})
to make meaningfull the operators $e^{itH}$; on the other hand
this implies also that the full Strichartz estimates \eqref{eq.strich}, 
\eqref{eq.strichF} hold for the group $e^{itH}$, hence we can write
\begin{equation}\label{eq.strV1}
    \|\Phi(v)\|_{\LII J p q}\leq C\;\|u_{0}\|_{L^{2}}+
       C\;\|W v\|_{\LII J {2}{\frac{2n}{n+2}}}+
       C\;\|F\|_{\LII J {\tilde p'}{\tilde q'}}
\end{equation}
for all admissible pairs $(p,q)$ and $(\tilde p,\tilde q)$. 
Notice that here $C$ is a constant depending on $V$ and 
the interval $J$ only, and can be assumed to be non 
increasing when $\delta\downarrow0$. This implies
\begin{equation}\label{eq.holdV}
    \|\Phi(v)\|_{\LII J p q}\leq C\;\|u_{0}\|_{L^{2}}+
       C\;\|W\|_{\LII J \infty {n/2}}  \|v\|_{\LII J {2} {\frac{2n}{n-2}}}+
       C\;\|F\|_{\LII J {\tilde p'}{\tilde q'}}
\end{equation}
and if $\delta$ is so small that
\begin{equation}\label{eq.condW}
    C\;\|W\|_{\LII J \infty {n/2}} \leq \frac12
\end{equation}
which is possible by the continuity of $V(t,x)$ as an 
$L^{n/2}$-valued function, we arrive at
\begin{equation}\label{eq.str5}
        \|\Phi(v)\|_{Z}\leq C\; \|u_{0}\|_{L^{2}}+
       \frac 1 2\|v\|_{Z}+
       C\; \|F\|_{\LII J {\tilde p'}{\tilde q'}}.
\end{equation}
This guarantees, as above, the existence of a unique local solution 
belonging to the space $Z$ and satisfying the Strichartz estimates 
with some constant $C(0)$ for some time interval $[0, \delta)$.

The same argument can be applied near each point $t_{0}\in I$. 
More precisely, let 
$J=[t_{0}-\delta,t_{0}+\delta]\cap I$ and assume $\delta>0$ 
is so small that the potential
$$W(t,x)=V(t,x)-V(t_{0},x)$$
satisfies
\begin{equation}\label{eq.t0}
    \|W(t,\cdot)\|_{L^{n/2}}\leq(2C(V(t_{0},x)))^{-1}\qquad\hbox{for }t\in J,
\end{equation}
where $C(V(t_{0},x))$ is the Strichartz constant corresponding 
to the potential $V(t_{0},x)$
and relative to the interval $[0,t_{0}+1]$.
Then we may argue as above, and we obtain that for any given 
initial time $t_{1}\in J$ and for any 
$f\in L^{2}$, the Cauchy problem
$$i\partial_{t}u-Hu=F(t,x)-W(t,x)u,\qquad u(t_{1})=f,
\qquad H=\Delta-V(t_{0},x)$$
(interpreted as usual in integral form via the group $e^{itH}$) 
has a unique solution
in $Z=C_{J}L^{2}\cap\LII J 2 {\frac{2n}{n-2}}$, which satisfies 
the Strichartz estimates
\begin{equation}\label{eq.str7}
        \|\Phi(v)\|_{Z}\leq 2C(t_{0})\; \|u_{0}\|_{L^{2}}+
            2C(t_{0})\; \|F\|_{\LII J {\tilde p'}{\tilde q'}}.
\end{equation}
for some constant $C(t_{0})$ depending on the point $t_{0}$ 
but \emph{not} on the initial time $t_{1}\in J$.

Now we may proceed by a continuation argument as follows. 
Extend the local solution constructed on $[0,\delta]$ to a 
maximal interval $[0,T^{*}[$; i.e., consider the union of 
all intervals $[0,\delta]$ on which a solution 
$u\in C([0,\delta];L^{2})\cap L^{2}(0,\delta;L^{\frac{2n}{n-2}})$ 
exists and satisfies (for all admissible pairs) the Strichartz 
estimates with some constant $C_{\delta}$. Assume by 
contradiction that $T^{*}<T$. Then the above local argument 
applied at $t_{0}=T^{*}$ on a suitable interval of the form 
$J=[T^{*}-\varepsilon,T^{*}+\varepsilon]$
shows that we can patch the maximal solution and extend 
it to $[0,T^{*}+\varepsilon]$. Moreover, we claim that the 
extended solution satisfies the Strichartz estimates on 
$[0,T^{*}+\varepsilon]$:
indeed, chosen any $t_{1}$ such that $T^{*}-\varepsilon<t_{1}<T^{*}$, 
by construction we see that 
the estimates hold both on $I_{1}=[0,t_{1}]$, with initial data at $t=0$:
\begin{equation}\label{eq.str8}
    \|u\|_{\LII {I_{1}} p q}\leq C'\;\|u(t_{0})\|_{L^{2}}
    +C'\;\|F\|_{\LII {I_{1}}{\tilde p'}{\tilde q'}},
\end{equation}
and on $J=[T^{*}-\varepsilon,T^{*}+\varepsilon]$, 
with initial data at $t=t_{1}$:
\begin{equation}\label{eq.str9}
    \|u\|_{\LII {J} p q}\leq C'\;\|u(t_{1})\|_{L^{2}}
    +C'\;\|F\|_{\LII J {\tilde p'}{\tilde q'}},
\end{equation}
for a suitable constant $C'$.
Since $\|u(t_{1})\|_{L^{2}}$ can be estimated exactly 
by \eqref{eq.str8} ($p=2,q=\infty$), we easily conclude 
the proof of our claim. This contradicts the assumption 
$T^{*}<T$ and we obtain that
$T^{*}=T$.

The modifications required to prove the final remark 
concerning the case $I=[0,\infty[$, and also Remark \ref{rem.nbd}, are obvious.

\section{Proof of the counterexamples}\label{sec.count}

\subsection{An eigenvalue problem.}\label{ssec.eigv} 
The first step in our construction requires to find a potential 
$V(x)$ such that
the operator $-\Delta+V(x)$ has a negative eigenvalue, i.e., 
such that the equation
\begin{equation}\label{eq.eigv}
    -\Delta u_{0}+ V(x)u_{0}+\gamma^{2}u_{0}=0
\end{equation}
admits a solution $u_{0}\in H^{1}$ for some $\gamma>0$.
There are many results on this problem, and in general there is 
a clear connection between
the number of such eigenvalues and the size of the 
negative part of $V$, in a suitable
norm. This is true both in the negative sense (explicit bounds 
on the number of the eigenvalues)
and in the positive sense, which is our main focus here. 
For instance, it is known that
(see \cite{ReedSimonIV}) if $V(x)\in L^{n/2}(\mathbb{R}^{n})$ 
satisfies the assumption
\begin{equation}\label{eq.neg}
    \hbox{the set
    $\{x\in\mathbb{R}^{n}\colon V(x)<0\}$ has a positive measure,}
\end{equation}
then there exists $\lambda_{0}>0$ such that, for all $\lambda>\lambda_{0}$,
the equation
\begin{equation}\label{eq.neg2}
    -\Delta u_{0}+ \lambda V(x)u_{0}+\gamma^{2} u_{0}=0
\end{equation}
admits at least a solution $f\in H^{1}$ for some $\gamma>0$. 
It can also be proved
that the dimension of the eigenspace grows to infinity as 
$\lambda$ tends to infinity.

However, for our purposes here we need only a much less 
precise result, which can be proved directly by an elementary 
variational argument. Indeed, take any smooth compactly 
supported function $w(x)$
such that $w(x_{0})>0$ at least in one point $x_{0}$. 
Then consider the minimization problem
with a constraint
\begin{equation}\label{eq.min}
    \min_{f \in M}\int_{\mathbb{R}^{n}}\left(|\nabla f|^{2}+|f|^{2}\right)\;dx
    \qquad\hbox{on}\qquad
    M=\left\{f\in H^{1}\colon \int_{\mathbb{R}^{n}}w(x)\;|f|^{2}dx=1\right\}.
\end{equation}
Note that $M$ is not empty, thanks to the assumption $w(x_{0})>0$. 
The existence of a solution to problem \eqref{eq.min} can be 
proved easily by a standard compactness argument, since we 
can work on the (bounded) support of $w(x)$. On the other hand,
the Euler-Lagrange equation of the problem is
\begin{equation}\label{eq.lagr}
    -\Delta f+f=\mu w(x)f
\end{equation}
(where $\mu$ is a Lagrange multiplier); hence, choosing 
$W(x)=-\mu w(x)$ and $u_{0}=f$, 
we see that $u_{0}$ solves the equation
\begin{equation}\label{eq.ell}
    -\Delta u_{0}+W(x)u_{0}+u_{0}=0
\end{equation}
and hence 
\begin{equation}\label{eq.ellt}
    u(t,x)=e^{-it}u_{0}(x)\hbox{\ \ solves\ \ } iu_{t}-\Delta u+W(x)u=0.
\end{equation}
Note also that a trivial bootstrapping argument gives 
$u_{0}\in H^{s}$ for all $s>0$.

This concludes the proof of Theorem \ref{th.3}, part (i).

\subsection{Proof of Theorem \ref{th.3}, case $1/r+n/(2s)<1$, 
$r\neq\infty$}\label{ssec.proof2a}

We start from the function \eqref{eq.ellt} and we apply the 
standard rescaling
\begin{equation}\label{eq.resc}
    u(t,x)\mapsto u_{\epsilon}(t,x)=u(\epsilon^{2}t,\epsilon x)\equiv
       e^{-i\epsilon^{2}t}u_{0}(\epsilon x).
\end{equation}
Then the function $u_{\epsilon}$ solves globally
\begin{equation}\label{eq.schep}
    i\partial_{t}u_{\epsilon}-\Delta u_{\epsilon}
    +W_{\epsilon}(x)u_{\epsilon}=0,\qquad
    W_{\epsilon}(x)=\epsilon^{2}W(\epsilon x).
\end{equation}
Consider now two monotone sequences of positive real numbers
\begin{equation}\label{eq.seq}
    0=T_{0}<T_{1}<\dots<T_{k}\uparrow+\infty, 
    \qquad 0<\epsilon_{k}\downarrow0,\qquad k=0,1,2,3,\dots
\end{equation}
and define a potential $V(t,x)$ on 
$[0,+\infty[\times \mathbb{R}^{n}$ by patching the potentials 
$V_{\epsilon}$ as follows:
\begin{equation}\label{eq.Vtot}
    V(t,x)=W_{\epsilon_{k}}(x)\qquad
    \hbox{for}\qquad t\in [T_{k},T_{k+1}[,\quad k=0,1,2,\dots.
\end{equation}
Thus $u_{\epsilon_{k}}$ solves the equation
\begin{equation}\label{eq.glob}
    i\partial_{t}u-\Delta u+V(t,x)u=0
\end{equation}
on the interval $[T_{k},T_{k+1}[$. 

Choose now $r$ and $s$ satisfying
\begin{equation}\label{eq.rs}
    \frac 1 r+\frac n {2s} <1,\qquad r\neq\infty,
\end{equation}
and assume we can choose the parameters $T_{k},\epsilon_{k}$ 
in such a way that 
\begin{equation}\label{eq.cond1}
    \|V\|_{L^{r}L^{s}}\leq \|W\|_{L^{s}}
    \sum_{k\geq0}(T_{k+1}-T_{k})^{1/r}\epsilon^{2-n/s}<\infty,
\end{equation}
then $V(t,x)\in L^{r}([0,+\infty[;L^{s})$.
On the other hand by Theorem \ref{th.1} we can extend 
(uniquely) $u_{\epsilon_{k}}$ to a global solution of 
\eqref{eq.glob} in $C([0,+\infty[;L^{2}( \mathbb{R}^{n}))$ 
which we shall denote by $u_{k}(t,x)$. Notice that, by the 
same theorem, we have
\begin{equation}\label{eq.L2const}
    \|u_{k}(t,\cdot)\|_{L^{2}}\equiv \hbox{const.}
    \equiv \|u_{\epsilon_{k}}(T_{k})\|_{L^{2}}\equiv
    \epsilon^{-n/2}\|u_{0}\|_{L^{2}}
\end{equation}
recalling the explicit expression \eqref{eq.resc} of 
$u_{\epsilon}$. On the other hand, we can write
\begin{equation}\label{eq.LpLq}
    \|u_{k}\|_{L^{p}(\mathbb{R}^{};L^{q})}\geq
    \|u_{k}\|_{L^{p}(T_{k},T_{k+1};L^{q})}\equiv
    \|u_{\epsilon_{k}}\|_{L^{p}(T_{k},T_{k+1};L^{q})}\equiv
    (T_{k}-T_{k+1})^{1/p}\epsilon_{k}^{-n/q}\|u_{0}\|_{L^{q}}
\end{equation}
by an elementary calculation. The Strichartz esimates are violated when
\begin{equation}\label{eq.unb}
    \frac { \|u_{k}\|_{L^{p}(\mathbb{R}^{};L^{q})}}
    {\|u_{k}(0)\|_{L^{2}}}\quad\hbox{is unbounded,}
\end{equation}
and this holds provided the parameters $T_{k},\epsilon_{k}$ 
satisfy the condition
\begin{equation}\label{eq.cond2}
     \frac { \|u_{k}\|_{L^{p}(\mathbb{R}^{};L^{q})}}
           {\|u_{k}(0)\|_{L^{2}}}\geq
     \frac { \|u_{\epsilon_{k}}\|_{L^{p}(T_{k},T_{k+1};L^{q})}}
     {\|u_{\epsilon_{k}}(0)\|_{L^{2}}}\equiv 
     \frac {\|u_{0}\|_{L^{q}}}{\|u_{0}\|_{L^{2}}}
     (T_{k}-T_{k+1})^{1/p}\epsilon_{k}^{\frac n 2 - \frac n q}   \to\infty.
\end{equation}
In conclusion, we only need to adjust the parameters 
\eqref{eq.seq} so to satisfy the two conditions 
\eqref{eq.cond1}, \eqref{eq.cond2}:
\begin{equation}\label{eq.cond3}
    \sum_{k\geq0}(T_{k+1}-T_{k})^{1/r}
    \epsilon_{k}^{2-n/s}<\infty,\qquad
    (T_{k}-T_{k+1})^{1/p}\epsilon_{k}^{\frac n 2 - \frac n q}   \to\infty,
\end{equation}
given an admissible pair $(p,q)$ with $p\neq\infty$ 
and $(r,s)$ as in \eqref{eq.rs}. With the special choices
\begin{equation}\label{eq.cho1}
    T_{0}=0,\quad T_{k+1}=T_{k}+k^{\alpha},\quad
    \epsilon_{0}=1,\quad \epsilon_{k}=k^{-\beta/2},
    \qquad k=1,2,3,\dots
\end{equation}
for some $\alpha,\beta>0$, the conditions reduce to
\begin{equation}\label{eq.cond4}
    \frac \alpha r +\beta\frac n {2s}<\beta-1,
    \qquad \frac \alpha p +\beta\frac n {2q}> \beta\frac n 4.
\end{equation}
Since $(p,q)$ is admissible, the second condition 
simplifies to $\alpha>\beta$, and rearranging the first one we are reduced to
\begin{equation}\label{eq.cond5}
    \frac {\alpha-\beta} r +\beta\left(\frac 1 r +
    \frac n {2s}\right)<\beta-1,\qquad
    \alpha>\beta.
\end{equation}
The term in brackets is smaller then 1 by assumption, 
hence if we choose any
\begin{equation}\label{eq.cond6}
    \alpha>\beta>\left[1-\left(\frac 1 r +\frac n {2s}\right)\right]^{-1}
\end{equation}
with $\alpha$ close enough to $\beta$,
we conclude the proof of the first part of Theorem \ref{th.2}, (ii).

\subsection{Proof of Theorem \ref{th.3}, case $1/r+n/(2s)>1$, 
$r\neq\infty$}\label{ssec.proof2b}

As in case \ref{ssec.proof2a} the proof is based on a 
rescaling argument. 
First of all we prove part (ii).
Consider again the rescaled solution \eqref{eq.resc} 
which solves equation \eqref{eq.schep} globally with a 
smooth compactly supported potential $W_{\epsilon}(x)=
\epsilon^{2}W(\epsilon x)$. Now, take two monotone 
sequences of positive real numbers
\begin{equation}\label{eq.seq2}
    1= \epsilon_{0}< \epsilon_{1}<\dots< \epsilon_{k}\uparrow+\infty, \qquad 
    0<\delta_{k}\downarrow0,\qquad k=0,1,2,3,\dots
\end{equation}
and define a potential $V(t,x)$ on $[0,+\infty[\times\mathbb{R}^{n}$ 
as follows:
\begin{equation}\label{eq.V0}
V(t,x) =
    \begin{cases}
       W_{\epsilon_{k}}(x) & \text{if $t\in[k,k+\delta_{k}]$, 
       $x\in\mathbb{R}^{n}$,}\\
        0 & \text{elsewhere.}
    \end{cases}
\end{equation}
Note that $V(t,x)\in \LI \infty\infty$ for any bounded time 
interval $I$, while globally
\begin{equation}\label{eq.cond7}
    \|V\|_{L^{r}L^{s}}\leq \|W\|_{L^{s}}
    \sum_{k\geq0}\delta_{k}^{1/r}\epsilon_{k}^{2-n/s}.
\end{equation}
As above, the function $u_{\epsilon_{k}}$ solves the equation
\begin{equation}\label{eq.glob2}
    i\partial_{t}u-\Delta u+V(t,x)u=0
\end{equation}
on the interval $t\in[k,k+\delta_{k}]$, and can be extended 
to a global solution $u_{k}(t,x)$ of the same equation thanks 
to the existence part of Theorem \ref{th.1} (recall that $V\in \LI\infty\infty$).
Moreover, $u_{k}$ has a conserved energy
\begin{equation}\label{eq.energy2}
    \|u_{k}(t,\cdot)\|_{L^{2}}\equiv\|u_{k}(k,\cdot)\|_{L^{2}}
    \equiv \epsilon^{-n/2}\|u_{0}\|_{L^{2}}.
\end{equation}
Then, as before, we can estimate
\begin{equation}\label{eq.cond8}
     \frac { \|u_{k}\|_{L^{p}(\mathbb{R}^{};L^{q})}}
        {\|u_{k}(0)\|_{L^{2}}}\geq
     \frac { \|u_{\epsilon_{k}}\|_{L^{p}(k,k+\delta_{k};L^{q})}}
     {\|u_{\epsilon_{k}}(0)\|_{L^{2}}}\equiv 
     \frac {\|u_{0}\|_{L^{q}}}{\|u_{0}\|_{L^{2}}}
     \delta_{k}^{1/p}\epsilon_{k}^{\frac n 2 - \frac n q} .
\end{equation}
Again, in order to violate the Stichartz estimates for an 
admissible couple $(p,q)$ and the potential $V\in L^{r}L^{s}$, 
it is sufficient to satisfy the two conditions
\begin{equation}\label{eq.cond9}
    \sum_{k\geq0}\delta_{k}^{1/r}\epsilon_{k}^{2-n/s}<\infty,\qquad
     \delta_{k}^{1/p}\epsilon_{k}^{\frac n 2 - \frac n q}\to\infty.
\end{equation}
With the special choices
\begin{equation}\label{eq.cho}
    \delta_{k}=k^{-\alpha},\qquad
    \epsilon_{k}=k^{\beta/2},
\end{equation}
the parameters $\alpha,\beta>0$ to be precised, we are reduced to
\begin{equation}\label{eq.cond10}
    -\frac\alpha r+\left(1-\frac n{2s}\right)\beta<-1,\qquad
    -\frac\alpha p+\left(\frac n 4-\frac n{2q}\right)\beta>0.
\end{equation}
Since $(p,q)$ is an admissible pair, the second condition is 
equivalent to $\alpha<\beta$ and we can rewrite the conditions as
\begin{equation}\label{eq.cond11}
    \frac{\alpha-\beta }r+\left(\frac 1 r +\frac n{2s}\right)\beta>\beta+1,\qquad
    \alpha<\beta.
\end{equation}
Recall now that we are considering the case
\begin{equation}\label{eq.caseb}
    \frac 1 r +\frac n{2s}>1,
\end{equation}
hence we may choose any $\beta$ such that
\begin{equation}\label{eq.condb}
    \beta>\left[\left(\frac 1 r +\frac n{2s}\right)-1\right]^{-1}
\end{equation}
and choosing then any $\alpha<\beta$ close enough to $\beta$, 
we easily satisfy the conditions \eqref{eq.cond11}.

This concludes the proof of part (ii) of Theorem \ref{th.3}.

Part (iii) can be proved by a simple 
modification of the preceding proof. Indeed, 
consider again the sequence $\delta_{k}=k^{-\alpha}$ 
constructed above, and notice that
it is not restrictive to assume that $\beta>\alpha>1$. 
Thus the series $\sum \delta_{k}$ converges, and the 
sequence of partial sums
\begin{equation}\label{eq.Tk}
    T_{k}=\sum_{j=0}^{k}\delta_{k}
\end{equation}
is positive, strictly increasing, and converges to
\begin{equation}\label{eq.T}
    \lim_{k\to\infty} T_{k}=T\equiv\sum_{k\geq0}\delta_{k}.
\end{equation}
We can now modify the definition 
\eqref{eq.V0} of the potential $V(t,x)$ as follows:
\begin{equation}\label{eq.V0bis}
V(t,x) =
    \begin{cases}
       W_{\epsilon_{k}}(x) & \text{if $t\in[T_{k},T_{k}+\delta_{k}]$, 
       $x\in\mathbb{R}^{n}$,}\\
        0 & \text{if $t\in[0,\delta_{0}[$.}
    \end{cases}
\end{equation}
This defines a potential on $I\times \mathbb{R}^{n}= 
[0,T]\times \mathbb{R}^{n}$, whose $\LI r s$ is 
given again by \eqref{eq.cond7}. The remaining arguments 
of the preceding proof apply without modification.

The proof of Theorem \ref{th.3} is concluded.

\subsection{Proof of Proposition \ref{prop.4}}\label{ssec.proof4}

The main tool of the proof is the \emph{pseudoconformal transform}
\begin{equation}\label{eq.conf}
    u(t,x)\quad\mapsto\quad
    U(T,X)=e^{-i\frac{|X|^{2}}{4T}}T^{-\frac n 2}u
    \left(-\frac 1 T,\frac X T\right)
\end{equation}
which takes a solution $u(t,x)$ of the Schr\"odinger 
equation in the variables $t,x$ into another solution 
of the same equation, in the variables $T,X$. If we apply 
the transform to the solution \eqref{eq.ellt}, we obtain a 
function $U(T,X)$ which solves
\begin{equation}\label{eq.schconf}
    i\partial_{T}U-\Delta_{X}U+V(T,X)U=0,\qquad U(1,X)=e^{i-i|X|^{2}}
    u_{0}(X),
\end{equation}
where the potential $V(T,X)$ is given by
\begin{equation}\label{eq.Vconf}
    V(T,X)=\frac1{T^{2}}W\left(\frac X T\right).
\end{equation}
It is easy to compute explicitly the norm of $V$ on the interval $[0,1]$:
\begin{equation}\label{eq.normV}
    \|V\|_{L^{r}(\delta,1;L^{s})}=\left(\int_{\delta}^{1} 
    T^{r(n/s-2)}dT\right)^{1/r}   \|W\|_{L^{s}}<\infty
\end{equation}
and this integral converges since our assumption 
\eqref{eq.condcontr} on the pair $(r,s)$ is equivalent to \begin{equation*}
    r\left(\frac n s -2\right)>-1.
\end{equation*}
On the other hand, the $\LI p q$ norm of $U(T,X)$ on 
an interval of the form $[\delta,1]$ with $0<\delta<1$ is given by
\begin{equation}\label{eq.normU}
    \|U\|_{\LI p q}=\left(\int_{\delta}^{1}
    T^{p\left(\frac n q -\frac n 2\right)}\right)^{1/p}\|W\|_{L^{q}}\equiv
    \left(\int_{\delta}^{1}T^{-2}\right)^{1/p}\|W\|_{L^{q}}
\end{equation}
since admissible pairs $(p,q)$ satisfy $p(n/q-n/2)\equiv-2$. 
This implies that $U(T,X)$ belongs ot all $\LI p q$ spaces for 
$I=[\delta,1]$ for all $0<\delta<1$, but not for $I=[0,1]$ 
where the integral diverges. Note also that
\begin{equation*}
    \|U(1,\cdot)\|_{L^{2}}\equiv\|u_{0}\|_{L^{2}}.
\end{equation*}

It is sufficient now to apply to $U(T,X)$ a reflection and a 
translation in time $T$ to obtain exactly the counterexample 
required in Theorem \ref{prop.4}. The proof is concluded.



\end{document}